\documentclass[letterpaper, 10 pt, conference]{ieeeconf}  

\IEEEoverridecommandlockouts                              

\usepackage{graphicx,amsmath}
\usepackage{amssymb}
\usepackage{amsfonts}
\usepackage{calrsfs}
\usepackage{mathrsfs}
\usepackage{caption}
\usepackage{subcaption}
\usepackage{multirow}
\usepackage{xcolor}
\usepackage{graphics}
\usepackage{float}
\usepackage{graphicx}
\usepackage{calrsfs}
\usepackage{mathtools}
\usepackage{amsmath}
\usepackage{psfrag}
\usepackage{blindtext}
 \newtheorem{teorema}{{\em Theorem}}{ }
 { }
 { }

\makeatletter
\renewcommand*\env@matrix[1][*\c@MaxMatrixCols c]{%
	\hskip -\arraycolsep
	\let\@ifnextchar\new@ifnextchar
	\array{#1}}
	\makeatother
\mathtoolsset{showonlyrefs}
 \usepackage{graphicx}
\usepackage{amssymb}
\usepackage{epstopdf}
\usepackage{xcolor,colortbl}
\usepackage{hyperref}
\usepackage{pbox}
\usepackage{epsfig,psfrag}
\usepackage{color}
\usepackage{amsmath}
\usepackage{mathrsfs}  
\usepackage{amsbsy}
\usepackage{lipsum} 
\usepackage{pgf}
\usepackage{tikz}
\usetikzlibrary{arrows, automata, shapes.multipart, chains, positioning, fit, graphs, spy, shapes, backgrounds, patterns,decorations.pathmorphing}
\usepackage{verbatim}
\tikzstyle{block} = [draw, fill=blue!20, rectangle, 
minimum height=1 cm, minimum width=2cm, rounded corners=0.1cm]
\tikzstyle{block2} = [draw, fill=black!15, rectangle, 
minimum height=0.75 cm, minimum width=1.5cm, rounded corners=0.1cm]
\tikzstyle{faultyblock} = [draw, fill=red!20, rectangle, 
minimum height=1 cm, minimum width=2cm, rounded corners=0.1cm]
\tikzstyle{reconblock} = [draw, fill=green!20, rectangle, 
minimum height=1 cm, minimum width=2cm, rounded corners=0.1cm]
\tikzstyle{sum} = [draw, fill=black!15, circle, minimum width=0.5cm]
\tikzstyle{input} = [coordinate]
\tikzstyle{output} = [coordinate]
\tikzstyle{pinstyle} = [pin edge={to-,thin,black}]
\usepackage{float}
\usepackage{isomath}
\usepackage{subcaption}

\title{\LARGE \bf
Multivariable Extremum Seeking Unit-Vector Control Design
}

\author{Enzo Ferreira Tomaz Silva$^{1}$, Pedro Henrique Silva Coutinho$^2$, Tiago Roux Oliveira$^2$, Miroslav Krsti\'{c}$^3$
\thanks{*This work was supported by the Brazilian agencies CNPq (Grant numbers: 407885/2023-4 and 309008/2022-0), CAPES, and FAPERJ.}
\thanks{$^{1}$ Enzo Silva is with the Graduate Program in Electronics Engineering, State University of Rio de Janeiro, Brazil.
			        {\tt\small enzotomazsilva@gmail.com} }
\thanks{$^{2}$ Pedro Coutinho and Tiago Roux Oliveira are with the Department of Electronics and Telecommunication
Engineering, State University of Rio de Janeiro, Brazil.
			        {\tt\small phcoutinho@eng.uerj.br, tiagoroux@uerj.br }}
\thanks{$^{3}$ M. Krsti\'{c} is with the  Department of Mechanical and Aerospace Engineering, University of California - San Diego, USA (e-mail: krstic@ucsd.edu).}
}

\begin{document}

\maketitle
\thispagestyle{empty}

\pagestyle{empty}

\begin{abstract}
This paper investigates multivariable extremum seeking using unit-vector control. By employing the gradient algorithm and a polytopic embedding of the unknown Hessian matrix, we establish sufficient conditions, expressed as linear matrix inequalities, for designing the unit-vector control gain that ensures finite-time stability of the origin of the average closed-loop error system. Notably, these conditions enable the design of non-diagonal control gains, which provide extra degrees of freedom to the solution. The convergence of the actual closed-loop system to a neighborhood of the unknown extremum point is rigorously proven through averaging analysis for systems with discontinuous right-hand sides. Numerical simulations illustrate the efficacy of the proposed extremum seeking control algorithm.
\end{abstract}

\section{Introduction}


Extremum seeking control is an adaptive, real-time, and model-free approach suitable to identify an optimal point where a given objective function, with unknown parameters, is either maximized or minimized, thereby reaching its extremum~\cite{ESsurvey,LivroTiagoES}.
A widely recognized extremum-seeking method is the gradient algorithm. This approach optimizes the system by introducing sinusoidal perturbations into the control scheme, allowing the algorithm to estimate the gradient direction and adjust the control signal accordingly to drive the system towards the desired optimum point.


This approach has been successfully employed to different classes of extremum-seeking control problems, such as extremum seeking control with delay systems~\cite{ArtDelay2}, maps in cascade with partial differential equations~\cite{ArtPDE, galvao2022extremum,silva2023extremum}, cooperative games with Nash equilibrium~\cite{rodrigues2024nash}, and networked control systems~\cite{rodrigues2025event}. However, the finite-time convergence, at least for the average dynamics, is not addressed in that work.

The unit-vector control (UVC) is a specific sliding mode control (SMC) approach used to ensure robust and finite-time convergence of the system's trajectories. The UVC is a nonlinear control method that modifies the dynamic structure of the system by adding a discontinuous control law that drives the system towards
the origin in finite time 
~\cite{SlidingMode2024,oliveira2023sliding}. This introduces interesting robustness properties for the closed-loop system. Due to these features, SMC-based approaches have been exploited in the context of extremum seeking control, such as for energy optimization~\cite{toloue2017multivariable}, Nash equilibrium seeking for quadratic duopoly game~\cite{rodrigues2024sliding}, extremum seeking using mixed integral sliding mode controller~\cite{solis2019extremum}, robust integral sliding mode control for optimization of measurable cost functions~\cite{solis2021robust}, source seeking for nonholonomic systems~\cite{mellucci2025source}, extremum-seeking for reaction systems with uncertainty estimation~\cite{lara2017robust}, switching-based extremum seeking approach~\cite{chen2017switching}, extremum seeking based on second order sliding modes~\cite{angulo2015nonlinear}, and multivariable extremum-seeking by periodic switching functions~\cite{peixoto2020multivariable}. 
More recently, fixed-time solutions for Nash Equilibrium seeking have been proposed by~\cite{poveda2022fixed} for time-invariant and non-smooth systems. Prescribed-time strategies (time-varying and smooth) were also proposed for extremum seeking~\cite{yilmaz2024prescribed} and source seeking~\cite{todorovski2023practical}.
Particularly, the control methods in~\cite{poveda2022fixed} exhibit similarities to those in higher-order SMC, such as the super-twisting algorithm~\cite{geromel2024multivariable}.
However, most approaches are based on the use of relay-type systems. Moreover, they do not provide constructive approaches
to designing the control gain, which is even more involved in the multivariable case.

Inspired by the discussion above, this paper deals with the design of a multivariable extremum seeking scheme using a unit-vector controller, rather than proportional control laws. This is the first effort to exploit 
a polytopic embedding of the Hessian matrix of the quadratic map and pursue a suitable transformation from which the average closed-loop dynamics can be rewritten in an appropriate form that enables deriving LMI-based conditions to design the control gain of the UVC law that renders finite-time stability.  Moreover, we also propose an optimization problem to incorporate an objective function related to the guaranteed reaching time minimization. 
Finally, by employing Lyapunov stability arguments and averaging theorem for systems with discontinuous right-hand sides~\cite{plotnikov1979averaging}, we guarantee the convergence to a neighborhood of the unknown extremum point. 

This paper is organized as follows. Section~\ref{sec:problem} provides the problem formulation, including the development of 
the closed-loop average dynamics and the polytopic embedding for the uncertain Hessian matrix.
Then, in Section~\ref{sec:main-results}, we provide a LMI-based condition to design the UVC gain,
and we employ averaging theorem arguments to assure the convergence of the extremum seeking control system
in the non-average sense. A numerical example is provided in Section~\ref{sec:results} to illustrate
the effectiveness of the proposed approach. Finally, Section~\ref{sec:conclusion} concludes this paper.

\textbf{Notation.}
$\mathbb{R}^n$ denotes the $n$-dimensional  Euclidean space and $\mathbb{R}^{m\times n}$ the set of real matrices $m\!\times\!n$. 
$X\!>\! 0~(X \!<\! 0)$ denotes $X$ is a symmetric positive (negative) definite matrix. 


\section{Problem Formulation}
\label{sec:problem}

\begin{figure}[!ht]
\begin{center}
\includegraphics[width=8.4cm]{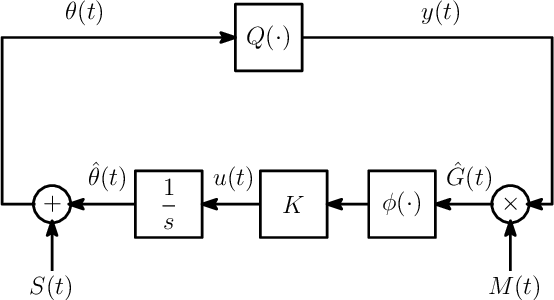}
\caption{Extremum seeking with unit vector control.} 
\label{fig:diag_ES_sliding}
\end{center}
\end{figure}

Consider the extremum seeking control with the unit vector control law based on the gradient algorithm shown in Figure~\ref{fig:diag_ES_sliding}.

We consider a quadratic static map given by
\begin{equation} \label{eq:map_static}
       y(t) = Q(\theta(t)) =  Q^* + \frac{1}{2}(\theta(t) - \theta^*)^\top H(\theta(t) - \theta^*),
\end{equation}
where $Q^* \in \mathbb{R}^n$ is the unknown optimal point of the map, 
$\theta^* \in \mathbb{R}^n$ is the unknown optimizer of the map, 
$\theta \in \mathbb{R}^n$ is the input vector, 
$H \in \mathbb{R}^{n \times n}$ is the unknown Hessian Matrix and 
$y \in \mathbb{R}^n$ is the map output.
Even though the Hessian matrix is unknown, it can be assumed that it is a positive definite matrix when the minimum point is desired.

In this scheme, the signal $\theta(t)$ that is applied to the static map is 
\begin{equation} \label{eq:input_map}
    \theta(t) = \hat{\theta}(t) + S(t)
\end{equation}
where $\hat{\theta}(t) \in \mathbb{R}^n$ is the estimated value of $\theta^*$, whose dynamics is described as follows
\begin{align}
    \dot{\hat{\theta}}(t) = K \phi(\hat{G}(t))
\end{align}
where $K$ is the control gain to be designed, $\phi(\cdot)$ is a nonlinear function defined by 
\begin{align}
    \phi(\hat{G}(t)) = \frac{\hat{G}(t)}{\|\hat{G}(t)\|},
\end{align}
and $\hat{G}(t)$ is the gradient estimate given by
\begin{equation}
    \hat{G}(t) = M(t)y(t).
\end{equation}
The perturbation signals of the extremum-seeking control system are defined as follows
\begin{equation} \label{eq:vetor_S}
    S(t) = \begin{bmatrix}
       a_1\sin{(\omega_1 t)} & \cdots & a_n\sin{(\omega_n t)} \\
    \end{bmatrix}^\top
\end{equation}
\begin{equation} \label{eq:vetor_M}
    M(t) = \begin{bmatrix}
       \tfrac{2}{a_1}\sin{(\omega_1 t)} & \cdots & \tfrac{2}{a_n}\sin{(\omega_n t)} \\
    \end{bmatrix}^\top
\end{equation}
where $a_i$, $i = 1,\ldots,n$, are positive scalars, and the frequencies of the perturbation signals are selected such that 
\begin{align} \label{eq:relacao_sinais_omega_pert}
    \omega_i = \omega_i' \omega, \quad i = 1,\ldots,n,
\end{align}
and $\omega_i' \notin \{\omega_j', \frac{1}{2}(\omega_j'+\omega_k'), \omega_k'\pm \omega_l' \}$,
for all $i,j,k = 1,\ldots,n$.

By defining the estimation error
\begin{equation} \label{eq:error_map}
    \tilde{\theta}(t) = \hat{\theta}(t) - \theta^*,
\end{equation}
its dynamics can be expressed as follows:
\begin{align}
    \dot{\tilde{\theta}}(t) = \dot{\hat{\theta}}(t) = \phi(\hat{G}(t)) = \frac{\hat{G}(t)}{\|\hat{G}(t)\|}.
\end{align}
Thus, it is clear that if the gradient estimate $\hat{G}(t)$ converges to zero possible, then
the estimation error $\tilde{\theta}(t)$ also converges to zero.

By substituting the expression for $y(t)$ in~\eqref{eq:map_static}, the 
gradient estimate signal can be expressed as:
\begin{equation}
    \hat{G}(t) = M(t)\left(Q^* + \frac{1}{2}(\tilde{\theta}(t) + S(t))^\top H(\tilde{\theta}(t) + S(t))\right),
\end{equation}
or still
\begin{multline} \label{eq:grad_simp}
     \hat{G}(t) = M(t)Q^* + \frac{1}{2}M(t)\tilde{\theta}^\top(t)H\tilde{\theta}(t)\\
     + M(t)S^\top(t)H\tilde{\theta}(t) + \frac{1}{2}M(t)S^\top(t)HS(t).
\end{multline}
By defining the matrix
\begin{equation} \label{eq:delta}
    \Omega(t) = M(t)S^\top(t)H,
\end{equation}
the multiplication in \eqref{eq:delta} results in a matrix of the following form:
\begin{equation} \label{eq:delta_final}
    \Omega(t) = H + \Delta(t)H,
\end{equation}
where
    $\Delta_{ii} = 1-\cos(2\omega_it)$,
    $\Delta_{ij} = \frac{a_{j}}{a_{i}}\cos(\omega_i-\omega_j) - \frac{a_{j}}{a_{i}}\cos(\omega_i+\omega_j)$.

Then, \eqref{eq:grad_simp} can be expressed as follows:
\begin{multline} \label{eq:grad_simp_2}
     \hat{G}(t) = M(t)Q^* + \frac{1}{2}M(t)\tilde{\theta}^\top H\tilde{\theta}(t) \\
     +  \Omega(t)\tilde{\theta}(t) + \frac{1}{2}\Omega(t)S(t).
\end{multline}
Finally, \eqref{eq:grad_simp_2} can still be rewritten as:
\begin{equation} \label{eq:grad_derivado}
     \dot{\hat{G}}(t) = Hu(t) + \Delta(t)Hu(t) + \varsigma(t)
\end{equation}
where
\begin{multline} \label{eq:resto_grad_derivado}
        \varsigma(t) = \dot{M}(t)Q^* + \dot{\Delta}(t)H\tilde{\theta}(t) + \frac{1}{2}H\dot{S}(t)+ \frac{1}{2} \dot{\Delta}(t)HS(t)  \\ + \frac{1}{2}\Delta(t)H\dot{S}(t) + [\frac{1}{2}M(t)\tilde{\theta}^\top H\tilde{\theta}(t)].
\end{multline}

\subsection{Time-Scale Change}

For the stability analysis of the closed-loop system, a change in the time scale is introduced here. 
From the relation of the disturbance signal frequencies \eqref{eq:relacao_sinais_omega_pert}, it is clear that there exists a common period $T$, given by
\begin{align}
    T = 2\pi \times \mathrm{LCM}\left\lbrace\frac{1}{\omega_i}\right\rbrace, \quad  i  = 1, 2, \ldots,n,
\end{align}
where $\mathrm{LCM}$ denotes the least common multiple. The change of time scale of the system in \eqref{eq:grad_derivado} consists of a transformation $\tau = \omega t$, where
\begin{align}
    \omega := \frac{2\pi}{T}.
\end{align}
Thus, the system \eqref{eq:grad_derivado} can be rewritten as
\begin{align} \label{eq:esc_temp_mudanca_Grad}
    \frac{d\hat{G}\left(\tau\right)}{d\tau} = \frac{1}{\omega}\mathcal{F}\left(\tau, \hat{G}, \tilde{\theta}, \frac{1}{\omega}\right) 
\end{align}
where 
\begin{align}\label{eq:esc_temp_mudanca_Grad_function}
  \mathcal{F}\left(\tau,  \hat{G}, \tilde{\theta}, \frac{1}{\omega}\right)   = Hu(\tau) + \Delta(\tau)Hu(\tau) + \varsigma(\tau)
\end{align}

\subsection{Average System}

By computing the average version of \eqref{eq:esc_temp_mudanca_Grad}, we have that:
\begin{align} \label{eq:Sistema_Medio_Derivado}
    \frac{d\hat{G}_\mathrm{av}(\tau)}{d\tau} = \frac{1}{\omega}\mathcal{F}_\mathrm{av}(\hat{G}_\mathrm{av}) 
\end{align}
where
\begin{align} \label{eq:Media_calculada_Sistema}
    \mathcal{F}_\mathrm{av}(\hat{G}_\mathrm{av}) = \frac{1}{T} \int_0^T   \mathcal{F}_\mathrm{av}(\delta, \hat{G}_\mathrm{av}, 0) d\delta.
\end{align}
The average of each term in~\eqref{eq:Media_calculada_Sistema} is given below:
\begin{align} 
    S_\mathrm{av}(\tau) &= \frac{1}{T} \int_0^T   S(\delta) d\delta = 0, \dot{S}_\mathrm{av}(\tau) = \frac{1}{T} \int_0^T   \dot{S}(\delta) d\delta = 0, \label{eq:Media_Sinal_S} \\
   M_\mathrm{av}(\tau) &= \frac{1}{T} \int_0^T   M(\delta) d\delta = 0, \dot{M}_\mathrm{av}(\tau) = \frac{1}{T} \int_0^T   \dot{M}(\delta) d\delta = 0, \label{eq:Media_Sinal_M} \\
    \Delta_\mathrm{av}(\tau) &= \frac{1}{T} \int_0^T   \Delta(\delta) d\delta = 0, \dot{\Delta}_\mathrm{av}(\tau) = \frac{1}{T} \int_0^T   \dot{\Delta}(\delta) d\delta = 0. \label{eq:Media_Sinal_Delta}
\end{align}
As a result, one can obtain
\begin{align*}
    \Omega_{\mathrm{av}}(\tau) = \frac{1}{T} \int_0^T \Omega(\delta) d\delta = H,  
    \dot{\Omega}_{\mathrm{av}}(\tau) = \frac{1}{T} \int_0^T \dot{\Omega}(\delta) d\delta = 0.
\end{align*}
Hence, the average system is finally given by
\begin{equation} \label{eq:grad_medio}
     \dot{\hat{G}}_{\mathrm{av}}(\tau) = \frac{1}{\omega} Hu_{\mathrm{av}}(\tau).
\end{equation}

By introducing the average nonlinear compensation
term $\phi_{\mathrm{av}}(\hat{G}_{\mathrm{av}}) = \hat{G}_{\mathrm{av}}/\|\hat{G}_{\mathrm{av}}\|$, the average closed-loop system can be expressed as 
\begin{align}\label{eq:average_system1}
    \dot{\hat{G}}_{\mathrm{av}}(\tau) = \frac{1}{\omega}H K \frac{\hat{G}_{\mathrm{av}}(\tau)}{\|\hat{G}_{\mathrm{av}}(\tau)\|}.
\end{align}
where the unit vector control law computed in terms of the average gradient estimate is
\begin{align} \label{eq:control_law_UVC}
    u_{\mathrm{av}}(\tau) =  K \frac{\hat{G}_{\mathrm{av}}(\tau)}{\|\hat{G}_{\mathrm{av}}(\tau)\|}.
\end{align}

In general, solutions available in the literature are developed for the stability analysis of extremum seeking control systems, assuming the knowledge of the sign of the Hessian matrix $H$. Based on this, a diagonal structure with the opposite sign is assigned to the gain matrix $K$. Although this approach requires little knowledge of the Hessian matrix $H$, it becomes difficult to design the gain matrix using constructive design conditions via LMIs. 

Instead, in this work, we assume that the Hessian matrix $H$ is unknown but takes values within a polytopic set $H \in \mathrm{co}\{H_i\}_{i = 1}^n$. Thus, it is possible to parameterize the Hessian matrix as follows:
\begin{align}\label{eq:polytpoic_H}
    H = H(\alpha) = \sum_{i=1}^N \alpha_i H_i,
\end{align}
where $\alpha = (\alpha_1,\ldots,\alpha_N)$ is the vector of uncertain parameters that belong to the unitary simplex
\begin{align}
    \Lambda = \left\lbrace \alpha \in \mathbb{R}^N : \sum_{i = 1}^N \alpha_i = 1, \; \alpha_i \geq 0, i = 1,\ldots,N \right\rbrace,
\end{align}
and $H_i \in \mathbb{R}^{n \times n}$, $i = 1,\ldots, N$ are the polytope vertices, that are known matrices. Thus, it is possible to write the average closed-loop system as:
\begin{align}\label{eq:dinamica_Gav_politopo}
    \dot{\hat{G}}_{\mathrm{av}}(\tau) = \frac{1}{\omega}H(\alpha) K \frac{\hat{G}_{\mathrm{av}}(\tau)}{\|\hat{G}_{\mathrm{av}}(\tau)\|},
\end{align}

The problem addressed in this paper is to design a control gain $K \in \mathbb{R}^{n \times n}$ such that the average closed-loop system~\eqref{eq:dinamica_Gav_politopo} is finite-time stable. Then, the averaging theorem is used to prove the closed-loop stability of system~\eqref{eq:grad_derivado}.


\section{Main Results}
\label{sec:main-results}

The main results of this paper are introduced in this section.
First, we propose an LMI-based condition to design the control gain of the
average extremum-seeking UVC system. Then, by employing the averaging theorem,
we ensure the stability of the extremum-seeking control system in the non-averaged sense.
Finally, we provide a convex optimization problem to design the control gain for a given
pre-specified reaching time for the average system.

\subsection{Controller Design Condition}

Let $z_{\mathrm{av}} = r(\hat{G}_{\mathrm{av}})\hat{G}_{\mathrm{av}}$ where $r(\hat{G}_{\mathrm{av}}) = 1/\sqrt{\|\hat{G}_{\mathrm{av}}\|}$. In $z_{\mathrm{av}}$-coordinates, the closed-loop system \eqref{eq:dinamica_Gav_politopo} can be rewritten as
\begin{align} \label{eq:L-coordinates}
    \dot{z}_{\mathrm{av}} = -\frac{1}{2}r(\hat{G}_{\mathrm{av}})\Pi_{\hat{G}_{\mathrm{av}}}  H(\alpha)Kr(\hat{G}_{\mathrm{av}})z_{\mathrm{av}} \\+ r(\hat{G}_{\mathrm{av}})  H(\alpha)Kr(\hat{G}_{\mathrm{av}})z_{\mathrm{av}}
\end{align}
where $\Pi_{\hat{G}_{\mathrm{av}}} = \hat{G}_{\mathrm{av}}\hat{G}^\top_{\mathrm{av}}/\|\hat{G}_{\mathrm{av}}\|^2$ is a projection matrix which satisﬁes the following properties: $\Pi_{\hat{G}_{\mathrm{av}}} = \Pi^\top_{\hat{G}_{\mathrm{av}}}, \Pi^2_{\hat{G}_{\mathrm{av}}} = \Pi_{\hat{G}_{\mathrm{av}}}$ and $\|\Pi_{\hat{G}_{\mathrm{av}}}\| = 1, \forall \hat{G}_{\mathrm{av}} \in \mathbb{R}^n$~\cite{geromel2024lmi}.

Consider the following Lyapunov function candidate~\cite{SlidingMode2024}:
\begin{align} \label{eq:lyapunov_L}
    V(z_{\mathrm{av}}) = z_{\mathrm{av}}^\top P z_{\mathrm{av}},
\end{align}
where $P = P^\top > 0$, which ensures that $V(z_{\mathrm{av}})$ is positive definite for all $z_{\mathrm{av}} \neq 0 \in \mathbb{R}^n$.
The theorem below provides a constructive condition based on LMIs for designing the control gain of the extremum-seeking control system.

\begin{teorema} \label{thm:2}
Given a scalar $\mu > 0$, if there exist symmetric matrices $X \in \mathbb{R}^{n \times n}$ and $M \in \mathbb{R}^{n\times n}$, and a full matrix $L \in \mathbb{R}^{m \times n}$, such that the following conditions hold:
\begin{align} \label{eq:first_condition_UVC}
    X > 0, \quad  M > 0,
\end{align}
\begin{equation} \label{eq:second_condition_UVC}
    \begin{bmatrix}
        H_iL + L^\top H_i^\top + \frac{\mu}{4}I + M & L^\top H_i^\top\\
        H_iL & \mu I
    \end{bmatrix} < 0, \forall i \in \mathbb{N} \leq N
\end{equation}
then, the origin of the closed-loop system \eqref{eq:dinamica_Gav_politopo} with $K = LX^{-1}$ converges in finite time.
\end{teorema}
\begin{proof}
    Assume that the conditions \eqref{eq:first_condition_UVC}-\eqref{eq:second_condition_UVC} hold.  From \eqref{eq:first_condition_UVC}, it follows that $X$ is a nonsingular matrix and there exists $X^{-1}$, since $X>0$. By multiplying the inequalities in \eqref{eq:second_condition_UVC} by $diag(X^{-1},I)$ on the left and its transpose on the right, it follows that
\begin{equation} \label{eq:proof_1_condition_UVC}
    \begin{bmatrix}
        PH_iK + K^\top H_i^\top P + \frac{\mu}{4}P^2 + Q & K^\top H_i^\top\\
        H_iK & \mu I
    \end{bmatrix} < 0, \forall i \in \mathbb{N} \leq N
\end{equation}
for all $i \in \mathbb{N} \leq N$, where $P = X^{-1}, K = LX^{-1}$, and $Q = X^{-1}MX^{-1}$. Since $B \in  \mathcal{co}\{B_i\}^N_{i=1} $, if multiplying \eqref{eq:proof_1_condition_UVC} by $\alpha_i$ and sum all the inequalities from $1$ to $N$, and  then applying Schur complement, it can be obtained 
\begin{equation} \label{eq:proof_2_condition_UVC}
    \frac{1}{\mu}K^\top H^\top HK + PHL + L^\top H^\top P + \frac{\mu}{4}P^2 + Q < 0
\end{equation}

Provided that
\begin{equation} \label{eq:proof_3_condition_UVC}
    -\frac{1}{2}K^\top H^\top \Pi_{\hat{G}_{\mathrm{av}}}  P - \frac{1}{2}P \Pi_{\hat{G}_{\mathrm{av}}} HK \leq \frac{1}{\mu}K^\top H^\top HK + \frac{\mu}{4}P^2,
\end{equation}
since
\begin{equation} \label{eq:proof_4_condition_UVC}
    (\frac{1}{\sqrt{\mu}}BK + \frac{\sqrt{\mu}}{2}\Pi_{\hat{G}_{\mathrm{av}}}  P)^\top (\frac{1}{\sqrt{\mu}}BK + \frac{\sqrt{\mu}}{2}\Pi_{\hat{G}_{\mathrm{av}}}  P) \geq 0
\end{equation}
and $\|\Pi_{\hat{G}_{\mathrm{av}}}\| = 1$, then it follows from \eqref{eq:proof_2_condition_UVC} and \eqref{eq:proof_3_condition_UVC} that
\begin{equation} \label{eq:proof_5_condition_UVC}
    -\frac{1}{2}K^\top H^\top \Pi_{\hat{G}_{\mathrm{av}}}  P - \frac{1}{2}P \Pi_{\hat{G}_{\mathrm{av}}} HK + S B K +K^\top B^\top P + Q < 0.
\end{equation}

By multiplying \eqref{eq:proof_5_condition_UVC} with $z_{\mathrm{av}}^\top r(\hat{G}_{\mathrm{av}})$ on the left and its transpose on the right, it follows that
\begin{align} \label{eq:proof_6_condition_UVC}
    \dot{V}(z_{\mathrm{av}}) < -z_{\mathrm{av}}^\top r(\hat{G}_{\mathrm{av}})Qr(\hat{G}_{\mathrm{av}})z_{\mathrm{av}} < 0,
\end{align}
with $V(z_{\mathrm{av}})$ deﬁned in \eqref{eq:lyapunov_L}. By following similar arguments of the proof of~\cite[Theorem~1]{SlidingMode2024}, it is possible to ensure that the origin is globally attractive. To show the finite-time convergence, notice that $z_{\mathrm{av}}^\top r(\hat{G}_{\mathrm{av}})Qr(\hat{G}_{\mathrm{av}})z_{\mathrm{av}} \geq \lambda_{min}(Q)\|z_{\mathrm{av}}\|^2/\|\hat{G}_{\mathrm{av}}\| = \lambda_{min}(Q)$, hence, it is possible to obtain from \eqref{eq:proof_6_condition_UVC} that the reaching time is upper bounded by
\begin{align} \label{eq:proof_7_condition_UVC}
    T_r \leq V_0/\lambda_{min}(Q),
\end{align}
where
\begin{align}
        V_0 = V(\hat{G}_{\mathrm{av}}(0)) = \hat{G}^\top_{\mathrm{av}}(0)P\hat{G}_{\mathrm{av}}(0)/\|\hat{G}_{\mathrm{av}}(0)\|, 
\end{align}
for all $\hat{G}_{\mathrm{av}}(0) \neq 0$. This concludes the proof.
\end{proof}

\subsection{Stability Analysis Using averaging theorem}

\begin{teorema} \label{thm:1}
    Take the average closed-loop dynamic of the gradient estimate subject to saturation~\eqref{eq:average_system1}.
    If the theorem conditions~\ref{thm:1}
    are satisfied, 
    so, for $\omega > 0$ sufficiently big, the equilibrium $\hat{G}_{\mathrm{av}} = 0$ converges in finite time and $\tilde{\theta}_{\mathrm{av}}(t)$ converges exponentially to zero. In particular, exist constants $m, M_\theta, M_y > 0$ such as
    \begin{align} \label{eq:theta_desigualdades}
        \|\theta(t) - \theta^\ast\| \leq \overline{\kappa} e^{-\eta t} + \mathcal{O}\left(a + \frac{1}{\omega}\right) \\
        \label{eq:y_desigualdades}
        |y(t) - Q^\ast| \leq  \kappa_y e^{-\eta t} + \mathcal{O}\left(a^2 + \frac{1}{\omega^2}\right),
    \end{align}
    where $a = \sqrt{\sum_{i=1}^n a_i^2}$, taking $a_i$ the defined constants in~\eqref{eq:vetor_S} and $\overline{\kappa}$ and $\kappa_y$ constants wich depends of the initial condition $\theta(0)$. 
\end{teorema}
\begin{proof}
    From the equation~\eqref{eq:grad_simp_2}, it can be obtained that
    \begin{align}
        \hat{G}_{\mathrm{av}}(\tau) = \frac{1}{\omega}H \tilde{\theta}_{\mathrm{av}}(\tau),
    \end{align}
    since the quadratic term $\frac{1}{2}M(t)\tilde{\theta}^\top H\tilde{\theta}(t)$ can be removed in a local analysis, and the other terms have zero average.
    
    Assuming the following Lyapunov function
    \begin{align}
        V(\tilde{\theta}) = \tilde{\theta}_{\mathrm{av}}^\top \overline{P} \tilde{\theta}_{\mathrm{av}},
    \end{align}
    where $\overline{P} = H^\top P H$ is a symmetric positive definite matrix, provided that $H$ and $P$ are symmetric and positive definite. Thus, it is possible to show that $\theta_{\mathrm{av}}(t)$ also converges in finite time to zero.
    As the differential equation in~\eqref{eq:esc_temp_mudanca_Grad} has discontinuity, due to the presence of the unit-vector function, \eqref{eq:esc_temp_mudanca_Grad_function} is $T$-periodic and Lipschitz continuous, it can be guaranteed from~\cite{plotnikov1979averaging} that
    \begin{align}
        \|\tilde{\theta}(t) - \tilde{\theta}_{\mathrm{av}}(t)\| \leq \mathcal{O}\left(\frac{1}{\omega}\right).
    \end{align}

    By applying the triangular inequality, it can be guaranteed that
    \begin{align} 
        \| \tilde{\theta}(t) \| 
        &\leq \overline{\kappa} e^{-\eta t} \| \tilde{\theta}_{\text{av}}(0) \| + \mathcal{O}\left(\frac{1}{\omega}\right).
    \end{align}

    Applying the average mean theorem~\cite{plotnikov1979averaging}, we can conclude that
    \begin{align}
        \| \hat{G}(t) - \hat{G}_{\text{av}}(t) \| \leq \mathcal{O}\left(\frac{1}{\omega}\right).
    \end{align}
    Analogously, applying the triangular inequality, we can obtain that
    \begin{align}
        \| \hat{G}(t) \| 
        &\leq \overline{\kappa} e^{-\eta t} \| \hat{G}_{\text{av}}(0) \| + \mathcal{O}\left(\frac{1}{\omega}\right).
    \end{align}
    From~\eqref{eq:input_map} and the definition of $\tilde{\theta}(t)$, we have
    \begin{align}
        \theta(t) - \theta^* = \tilde{\theta}(t) + S(t).
    \end{align}
    Based on that, the following relation can be obtained:
    \begin{align} \label{eq:demonstracao_estabilidade_theta}
        \| \theta(t) - \theta^* \| \leq (\overline{\kappa}) e^{-\eta t} \|\theta(0)- \theta^*\| + \mathcal{O}\left(a + \frac{1}{\omega}\right)
    \end{align}

By defining the error variable
\begin{align}
    \tilde{y}(t) := y(t) - Q^*,~~~~y(t) = Q(\theta(t)),
\end{align}
calculating its norm, and using the Cauchy–Schwarz's inequality, gets
\begin{align}
    |\tilde{y}(t)| & = |y(t) - Q^*| = |(\theta(t) - \theta^*)^\top H(\theta(t) - \theta^*)| \\ &\leq \|H\| \|((\theta(t) - \theta^*))\|^2.
\end{align}
From~\eqref{eq:demonstracao_estabilidade_theta}, it is still possible to obtain
\begin{align}
    |\tilde{y}(t)| \leq \|H\| ((\overline{\kappa})^2 e^{-2\eta t} \|\theta(0)- \theta^*\|^2 + \mathcal{O}\left(a^2 + \frac{2a}{\omega} + \frac{1}{\omega^2}\right)
\end{align}

As $e^{-\eta t} \geq e^{-2 \eta t}$ for $\omega > 0$, and $a^2+\frac{1}{\omega^2} \geq \frac{2a}{\omega}$, by the Young's inequality, obtains
\begin{align}
    |y(t) - Q^*| \leq \kappa_y e^{-\eta t} + \mathcal{O}\left(a^2 + \frac{1}{\omega^2}\right),
\end{align}
where 
\begin{align*}
    \kappa_y = \|H\| (\overline{\kappa})^2  \|\theta(0)- \theta^*\|^2
\end{align*}

As a result, the inequalities \eqref{eq:theta_desigualdades} and \eqref{eq:y_desigualdades} are guaranteed. This concludes the proof.
\end{proof}

\subsection{Optimization problem to minimize the reaching time}

This section follows similar steps as~\cite[Section 3.2]{SlidingMode2024} to formulate an 
optimization for designing the UVC gain by minimizing the reaching time to the average gradient estimate convergence.
We first introduce the following constraint:
\begin{align} \label{eq:proof_1_condition2_UVC}
    \begin{bmatrix}
        \varphi & I \\
        I & X
    \end{bmatrix} \geq 0.
\end{align}
The condition in \eqref{eq:proof_1_condition2_UVC} implies from Schur complement that $\varphi \geq X^{-1}$. Thus, $P \leq \varphi$ or still $V(z_{\mathrm{av}}) \leq \varphi z_{\mathrm{av}}^\top z_{\mathrm{av}}$. Therefore, one can conclude that $\mathcal{B} \subset \mathcal{V}$, where 
\begin{align}\label{eq:conjunto_omega}
        \mathcal{V} = \left\lbrace \hat{G}_{\mathrm{av}} \in \mathbb{R}^n :  V(\hat{G}_{\mathrm{av}}) \leq 1 \right\rbrace,
\end{align}
Thus, if $\varphi$ is minimized, the set $\mathcal{V}$ is increased. 

For a given initial condition $\hat{G}_{\mathrm{av}}(0)$ (associated to $V_0$), 
the reaching time can be minimized by maximizing the
smallest eigenvalue of $Q$. This objective can be achieved with the constraint
\begin{align}\label{eq:maximize_Q}
\begin{bmatrix}
    M & X \\
    X & \rho I
\end{bmatrix} \geq 0.
\end{align}
From~\eqref{eq:maximize_Q}, it follows from Schur complement that
$M - \rho^{-1} X^2 \geq 0$. By multiplying both sides by $X^{-1}$,
we have that $Q \geq \rho^{-1} I$, since $Q = X^{-1} M X^{-1}$.
Thus, by minimizing $\rho$, the eigenvalues of $Q$ are maximized, thus
reducing the reaching time $T_r$.
Therefore, if $\hat{G}_{\mathrm{av}}(0)$ is taken inside of $\mathcal{V}$, the reaching time for the average dynamics is constrained by
\begin{align} \label{eq:reaching_time_bounded}
    T_r \leq V_0/\lambda_{min}(Q) \leq \rho
\end{align}
The optimization problem for minimizing the estimated reaching time for a given set of initial conditions, where $\varphi > 0$, is formulated as follows: 
\begin{align} \label{eq:optmization_problem}
    &\min~~\rho \\
    &\mathrm{subject~to~LMIs~in}~\eqref{eq:first_condition_UVC},~\eqref{eq:second_condition_UVC},~\eqref{eq:proof_1_condition2_UVC},\eqref{eq:maximize_Q}.
\end{align}

\section{Numerical Results}
\label{sec:results}
Consider the extremum-seeking control system with non-linear map \eqref{eq:map_static} with unknown Hessian matrix taking values in the polytopic set given by the following vertices

\begin{align}
    H_1 = (1-\overline{\delta}) H_0, \quad  H_2 = (1+\overline{\delta}) H_0,
\end{align}
where $\overline{\delta} > 0$ is a parameter used to construct the vertices of the polytopic and $H_0$ is the Hessian matrix used in~\cite{MultVar_ESC}:
\begin{align}
    H_0 = \begin{bmatrix}
        100 & 30 \\
        30 & 20
    \end{bmatrix} > 0.
\end{align}
In addition, for the simulation, it was assumed that unkown optimical points are $Q^* = 10$ and
$\theta^*$ =
$\begin{bmatrix}
    2 & 4 \\
\end{bmatrix}^\top$. 
Note that the unknown parameters $Q^*$ and $\theta^*$ are not used to design the control gain, only for system simulation.

The project was performed solving the optimzation problem \eqref{eq:optmization_problem},
considering $\overline{\delta} = 0.1$, $\varphi = 0.4$ and $\mu = 32.9034$. The controller designed was
$$K = \begin{bmatrix}
    -0.2393  &  0.3589 \\
    0.3589 &  -1.1965
\end{bmatrix}.$$

\begin{figure}[ht!]
     \centering
    \subfloat[\label{fig:sinal_controle_uvc}$u(t)$ -- {Theorem~\ref{thm:2}}]{
    \includegraphics[width=0.7\columnwidth]{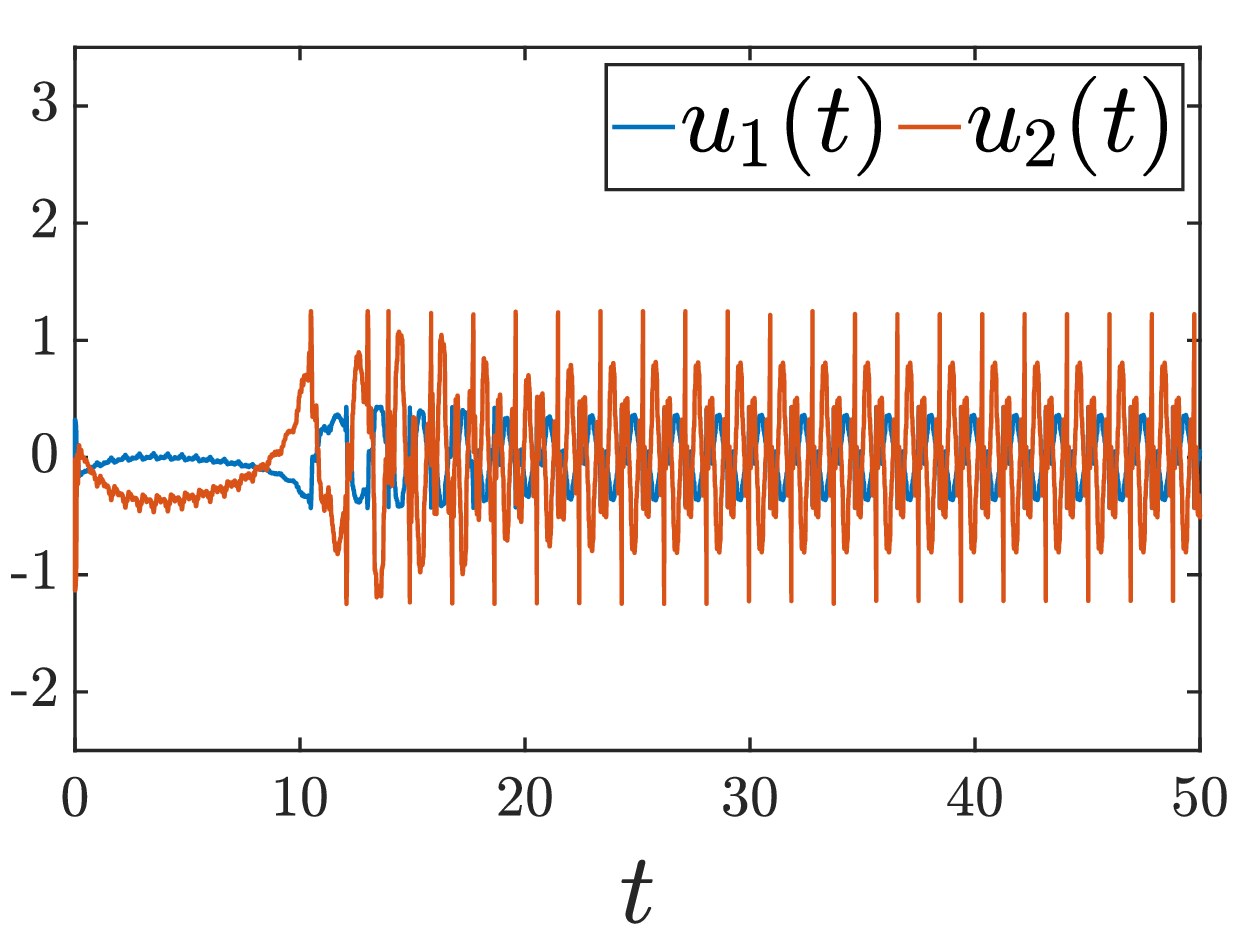}}

    \subfloat[\label{fig:sinal_theta_uvc}$\theta(t)$ -- {Theorem~\ref{thm:2}}]{
    \includegraphics[width=0.7\columnwidth]{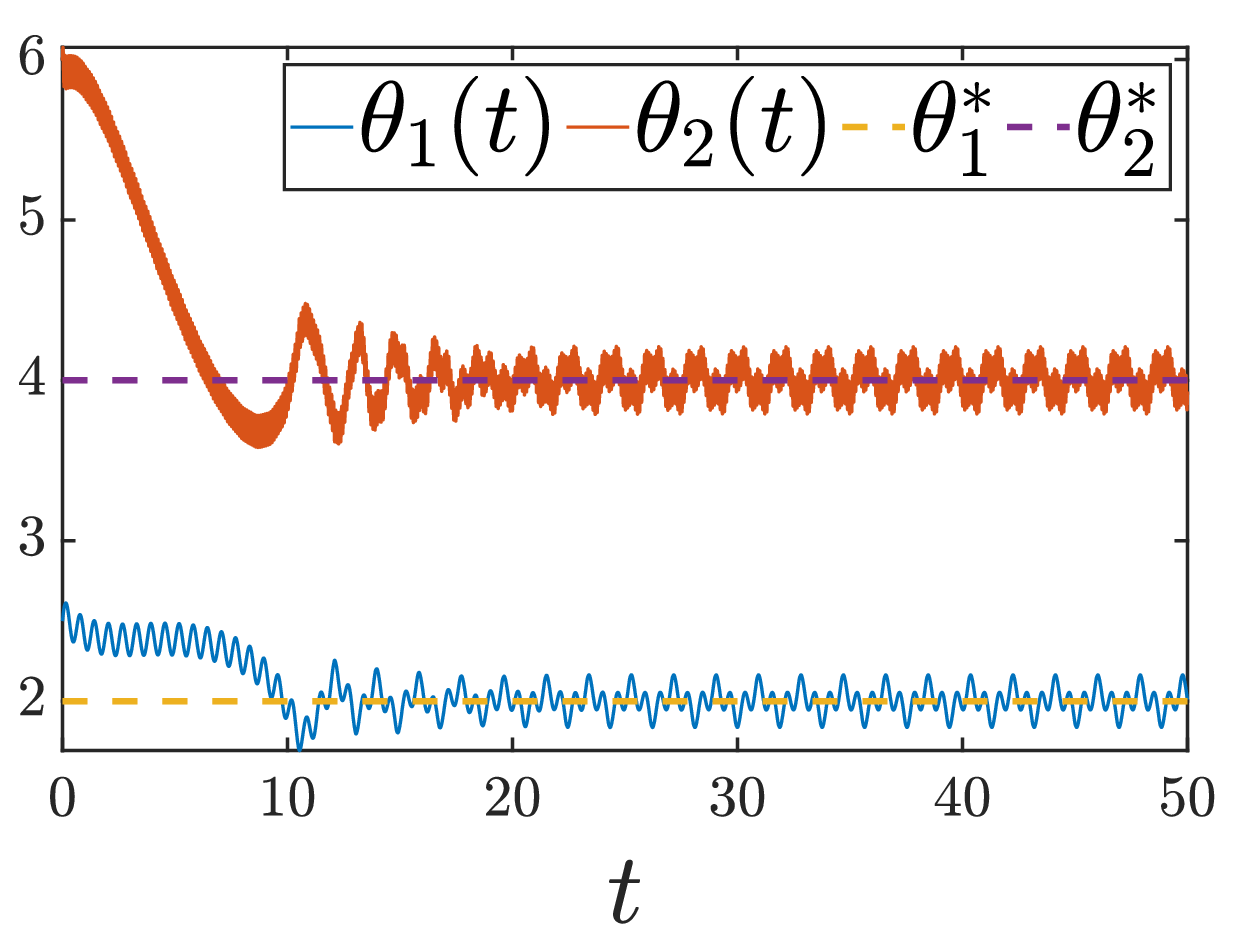}}

    \subfloat[\label{fig:sinal_saida_uvc}$y(t)$ -- {Theorem~\ref{thm:2}}]{
    \includegraphics[width=0.7\columnwidth]{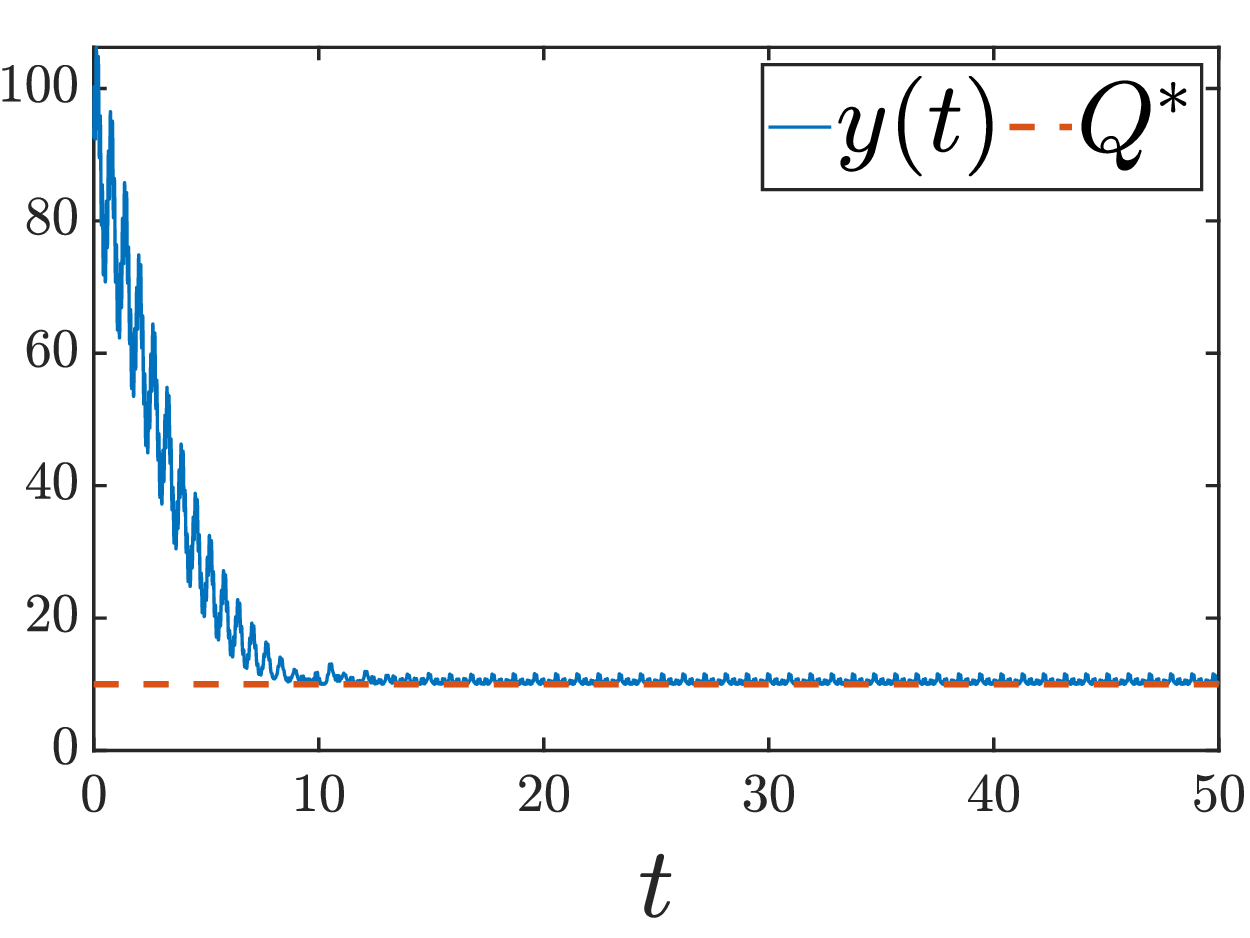}}
\caption{Response of the closed-loop system with the UVC law designed with Theorem~\ref{thm:2}.}
\label{fig:comparacao_ganho_LMI_arbitrario_Diag}
\end{figure}

For the simulation, the frequencies of the perturbation vectors~\eqref{eq:vetor_S} and~\eqref{eq:vetor_M} are
selected as $\omega_1 = 10$ rad/s and $\omega_2 = 70~\mathrm{rad/s}$, and their amplitudes are $a_1=a_2=0.1$.
The extremum-seeking control system was simulated with the designed controller using the Theorem conditions~\ref{thm:2}, developed in this work. The closed-loop simulation is depicted
in Figure~\ref{fig:comparacao_ganho_LMI_arbitrario_Diag}. 
The simulations were performed considering the initial condition $\theta(0) =  [2.5 \; 6]^\top$. 
As the results indicate, the system converges towards the optimum point using the designed UVC gain.
Thus, it was possible to ensure the extremum-seeking control system convergence with the proposed approach. Note also that the simulation was carried out considering a randomly generated value of $H$ taken inside the polytopic domain, which illustrates the robustness of the designed controller, as expected by the proposed
robust control design condition.

\section{Conclusion}
\label{sec:conclusion}
This paper tackles the problem of multivariable extremum control using a unit-vector control law.  By assuming an uncertain polytopic representation of the Hessian matrix, a constructive LMI-based condition is derived for designing the unit vector control gain that ensures the finite-time convergence of the average system,
instead of exponential stability, usually pursued in classical extremum
seeking with proportional control laws.
Then, by applying the averaging theorem for systems with discontinuous right-hand sides, it is shown that the system trajectories converge to a region around the unknown optimal point.
One of the main contributions of this conference paper is a constructive approach to designing the control gain of the extremum-seeking control system, which can be extended to derive novel conditions that ensure
stronger convergence guarantees, such as fixed-time convergence as in~\cite{poveda2022fixed}. Future researches lie in the design and analysis of different control problems with multivariable unit-vector approach, as considered in the following references \cite{paper1,paper2,paper3,paper4,paper5,paper6,paper7,paper8,paper9,paper10,paper11,paper12,paper13,paper14,paper15,paper16,paper17,paper18,paper19,paper20}.


\bibliographystyle{IEEEtran} 
\bibliography{references}


\end{document}